\documentclass[11pt]{article}
\usepackage{amstext,amssymb,amsmath,amsbsy}

\usepackage{hyperref}
\usepackage{amscd}
\usepackage{amsfonts}
\usepackage{indentfirst}
\usepackage{verbatim}
\usepackage{amsmath}
\usepackage{amsthm}
\usepackage{enumerate}
\usepackage{graphicx}
\usepackage{color}
\usepackage[OT1]{fontenc}
\usepackage[latin1]{inputenc}
\usepackage[english]{babel}
\usepackage{amssymb}

\newtheorem{theorem}{Theorem}
\newtheorem{lemma}{Lemma}

\newtheorem{remark}{Remark}

\numberwithin{equation}{section}

\newcommand{\proofend}{\hfill $\Box$ }

\newcommand{\dsp}{\displaystyle}

\newcommand{\epss}{s}
\newcommand{\eps}{\varepsilon}
\newcommand{\dive}{\operatorname{div}}

\newcommand{\bH}{{\bf H}}

\newcommand{\Div}{{\rm div}}
\newcommand{\ii}{i}
\newcommand{\pB}{\partial B}

\newcommand{\beq}{\begin{equation}}
\newcommand{\eeq}[1]{\label{#1}\end{equation}}
\newcommand{\Beq}{\begin{equation}\left\{\begin{array}{rcll}}
\newcommand{\Eeq}[1]{\end{array}\right.\label{#1}\end{equation}}
\newcommand{\Eeqq}{\end{array}\right.\nonumber\end{equation}}

\newcommand{\ds}{\displaystyle}
\newcommand{\idelta}{\ii \delta}

\newcommand{\mR}{\mathbb{R}}
\newcommand{\mZ}{\mathbb{Z}}
\newcommand{\mC}{\mathbb{C}}

\newcommand{\testfn}{e^{-\beta r^{-p}}}
\newcommand{\testfnn}{e^{2\beta r^{-p}}}
\newcommand{\X}{{B_{R_3}\setminus  \overline B_{R_1}}}
\newcommand{\XX}{\overline B_{R_3} \setminus B_{R_1}}

\newcommand{\pX}{\partial (\X)}
\newcommand{\annulus}[2]{B_{{#1}} \setminus \overline B_{{#2}}}

\newcommand{\supp}{\mbox{supp }}

\numberwithin{equation}{section}

\title{Cloaking using complementary media for the Helmholtz equation and a three spheres inequality for second order elliptic equations}
\author{Hoai-Minh Nguyen   \footnote{Mathematics Section, \'Ecole Polytechnique F\'ed\'erale de Lausanne, Station 8,  CH-1015 Lausanne, Switzerland, hoai-minh.nguyen@epfl.ch} \footnote{The research is partially supported by NSF grant DMS-1201370 and by the Alfred P. Sloan Foundation.} \, and Loc Hoang Nguyen  \footnote{Mathematics Section, \'Ecole Polytechnique F\'ed\'erale de Lausanne, Station 8,  CH-1015 Lausanne, Switzerland, loc.nguyen@epfl.ch}}

\begin{document}

\maketitle

\begin{abstract}
Cloaking using complementary media was suggested by Lai et al. in \cite{LaiChenZhangChanComplementary}. This was proved in \cite{Ng-Negative-Cloaking} in  the quasistatic regime.  One of the difficulties in the study of this problem is the appearance of the localized resonance, i.e., the fields  blow up in some regions and remain bounded in some others as the loss goes to 0. To this end,  the author introduced the technique of removing localized singularity and used a standard three spheres inequality.  The method used in \cite{Ng-Negative-Cloaking} also works for the Helmholtz equation. However, it requires small size of the cloaked region for large frequency due to the use of the (standard) three spheres inequality.  In this paper, we give a proof of cloaking using complementary media in the finite frequency regime without imposing any condition on  the cloaked region; hence the cloak works for all frequency. To successfully apply the approach in \cite{Ng-Negative-Cloaking}, we establish a new three spheres inequality. A modification of the cloaking setting to obtain illusion optics is also discussed. 
\end{abstract}



\section{Introduction}\label{sec: intro}

Negative index materials (NIMs) were investigated theoretically by Veselago in \cite{Veselago}. The  existence of such materials was confirmed by Shelby, Smith, and Schultz in \cite{ShelbySmithSchultz}. 
The study of NIMs has attracted a lot attention in the scientific community thanks to their interesting properties and  applications. One of the appealing one is cloaking using complementary media.

Cloaking using NIMs or more precisely cloaking using complementary media was suggested by Lai et al in  \cite{LaiChenZhangChanComplementary}. Their work was inspired from the notion of complementary media suggested by  Pendry and Ramakrishna  in \cite{PendryRamakrishna0}.  Cloaking using complementary media was established in \cite{Ng-Negative-Cloaking} in the quasistatic regime using  slightly different schemes from \cite{LaiChenZhangChanComplementary}.  Two difficulties in the study of cloaking using complementary media are as follows. Firstly, this problem is unstable since the equations describing the phenomenon  have sign changing coefficients, hence the ellipticity is lost. Secondly,  the localized resonance, i.e., the field blows up in some regions and remains bounded in some others,  might appear.  To handle these difficulties, in \cite{Ng-Negative-Cloaking} the author introduced the technique of removing localized singularity and used a standard three spheres inequality. The approach in \cite{Ng-Negative-Cloaking} also involved the reflecting technique introduced in \cite{Ng-Complementary}.  The method in \cite{Ng-Negative-Cloaking} also works for the Helmholtz equation; however since the largest radius in the (standard) three spheres inequality is  small as frequency is large  (see Section~\ref{Sect-ThreeSphereInequality} for further discussion), the size of the cloaked region is required to be small for large frequency.

In this paper, we present a proof of cloaking using complementary media in the finite frequency regime. Our goal is not to  impose any condition on the size of the cloaked region (Theorem~\ref{thm1}); hence the cloak works for all frequency. To successfully apply the approach in \cite{Ng-Negative-Cloaking}, we establish a new three spheres inequality for the second order elliptic equations which holds for arbitrary radius (Theorem~\ref{thm-threesphere} in Section~\ref{Sect-ThreeSphereInequality}). This inequality is inspired from the unique continuation principle and its proof is in the spirit of  Protter  in \cite{Protter60}.  A modification of the cloaking setting to obtain illusion optics  is  discussed  in  in Section~\ref{sect-illusion} (Theorem~\ref{thm2}).  This involves the idea of superlensing in \cite{Ng-Superlensing}. Cloaking using complementary media for electromagnetic waves is investigated in \cite{Ng-Negative-Cloaking-Maxwell}.

In addition to cloaking using complementary media, other application of NIMs are superlensing using complementary media as suggested in \cite{PendryNegative, PendryCylindricalLenses, PendryRamakrishna} (see also \cite{NicoroviciMcPhedranMilton94}) and confirmed in \cite{Ng-Superlensing, Ng-Superlensing-Maxwell},  and cloaking via anomalous localized resonance \cite{MiltonNicorovici} (see also \cite{AmmariCiraoloKangLeeMilton, KohnLu, Ng-CALR}). Complementary media 
were studied in a general setting in \cite{Ng-Complementary, Ng-Superlensing-Maxwell} and  play an important role in these applications see \cite{Ng-CALR-CRAS, Ng-Superlensing, Ng-CALR, Ng-Negative-Cloaking, Ng-Superlensing-Maxwell, MinhLoc}.

Let us describe the problem more precisely. Assume that the cloaked region is the annulus  $B_{\gamma r_2} \setminus B_{r_2}$ for some $r_{2}> 0$ and $1 < \gamma < 2$ in which the medium is characterized by a matrix $a$ and a function $\sigma$. The assumption on the cloaked region by all means imposes no restriction  since any bounded set is a subset of such a region provided that the radius and the origin are appropriately chosen.  The idea suggested by Lai et al. in \cite{LaiChenZhangChanComplementary} in two dimensions is to construct a complementary media in $B_{r_2} \setminus B_{r_1}$ for some $0 < r_1 < r_2$.

 In this paper, instead of taking the schemes of Lai et al., we use a scheme from \cite{Ng-Negative-Cloaking} which is inspired but different from the ones from \cite{LaiChenZhangChanComplementary}. Following \cite{Ng-Negative-Cloaking}, the cloak contains two parts. The first one, in $B_{r_2} \setminus  B_{r_1}$, makes use of complementary media to cancel the effect of the cloaked region and the second one, in $B_{r_1}$, is to fill the space which ``disappears" from the cancellation by the homogeneous media. Concerning the first part, 
instead of $B_{\gamma r_2} \setminus B_{r_2}$, we consider $B_{r_3} \setminus B_{r_2}$ with $r_3 =  2 r_2$ (the constant $2$ considered here is just a matter of simple representation) as the cloaked region in which the medium is given by 
\begin{equation*}
\hat a, \hat \sigma = \left\{ \begin{array}{cl} a, \sigma & \mbox{ in } B_{\gamma r_2} \setminus B_{r_2}, \\[6pt]
I, 1 & \mbox{ in } B_{r_3} \setminus B_{\gamma r_2}. 
\end{array} \right. 
\end{equation*} 
The complementary media in $B_{r_2} \setminus B_{r_1}$ is given by 
\begin{equation*}
- F^{-1}_*\hat a, -  F^{-1}_*\hat \sigma, 
\end{equation*}
where $F: B_{r_2} \setminus \bar B_{r_1} \to B_{r_3} \setminus \bar B_{r_2}$ is the Kelvin transform with respect to  $\partial B_{r_2}$, i.e., 
\begin{equation}\label{def-F}
F(x) = \frac{r_2^2}{|x|^2} x. 
\end{equation}
Here
\begin{equation*}
T_*\hat a (y) = \frac{D T  (x)  \hat a(x) D T ^{T}(x)}{J(x)} \quad \mbox{ and } \quad T_* \hat \sigma (y) = \frac{\hat \sigma(x)}{J(x)},
\end{equation*}
where $x =T^{-1}(y)$ and $J(x) = |\det D T(x)|$  for a diffeomorphism $T$. It follows that 
\begin{equation}\label{cond-r}
 r_1 = r_2^2/ r_3.   
\end{equation}
Concerning the second part, the medium in $B_{r_1}$ is given by 
\begin{equation}\label{choice-A}
\Big(r_3^2/r_2^2 \Big)^{d-2} I, \Big(r_3^2/r_2^2 \Big)^{d}. 
\end{equation}
The reason for this choice will be explained later.

With the loss, the medium is characterized by $s_\delta A, s_0 \Sigma$, 
where 
\begin{equation}\label{defA}
A, \Sigma = \left\{ \begin{array}{cl} 
\hat a, \hat \sigma & \mbox{ in } B_{r_3} \setminus B_{r_2}, \\[6pt]
F^{-1}_*\hat a, F^{-1}_* \hat  \sigma & \mbox{ in } B_{r_2} \setminus B_{r_1}, \\[6pt]
\Big(r_3^2/r_2^2 \Big)^{d-2} I, \Big(r_3^2/r_2^2 \Big)^{d} & \mbox{ in } B_{r_1},\\[6pt]
I, 1 & \mbox{ otherwise}, 
\end{array} \right. 
\end{equation} 
and 
\begin{equation}\label{def-ss}
s_\delta = \left\{ \begin{array}{cl} -1 + i \delta  & \mbox{ in } B_{r_2} \setminus B_{r_1}, \\[6pt]
1 & \mbox{ otherwise}. 
\end{array} \right. 
\end{equation}
Physically, the imaginary part of $s_\delta A$ is the loss of the medium (more precisely the loss of  the medium in 
$B_{r_2} \setminus B_{r_1}$). Here and in what follows, we assume that,
\begin{equation}\label{cond-a}
\frac{1}{\Lambda} |\xi|^{2} \le a (x) \xi \cdot \xi  \le \Lambda |\xi|^{2} \quad \forall \, \xi \in \mR^{d}, \mbox{ for a.e. } x \in B_{\gamma r_{2}} \setminus B_{r_2}, 
\end{equation}
for some $\Lambda \ge 1$. In what follows, we assume in addition that 
\begin{equation}
	\hat a \, \mbox{is Lipschitz in}\, B_{r_3} \setminus B_{r_1}.
\end{equation}
We can verify that  medium $s_0 A$ is of reflecting complementary property, a concept introduced in \cite[Definition 1]{Ng-Complementary}, by considering  diffeomorphism  $G: \mR^{d} \setminus \bar B_{r_{3}} \to B_{r_{3}} \setminus \{0\}$ which is the Kelvin transform with respect to $\partial B_{r_{3}}$, i.e., 
\begin{equation}\label{def-G}
G(x) = r_{3}^{2} x/ |x|^{2}. 
\end{equation}
It is important to note that  
\begin{equation}\label{GF}
G_{*} F_{*} A = I \mbox{ in } B_{r_{3}}
\end{equation}
since $G\circ F (x) = (r_3^2/ r_2^2)x$. This is the reason for choosing $A$ in  \eqref{choice-A}.  
  
Let $\Omega$ be a smooth open subset of $\mR^d$ ($d=2, \, 3$) such that $B_{r_3} \subset \subset \Omega$. Given $f \in L^{2}(\Omega)$, let  $u_\delta, \, u \in H^1_{0}(\Omega)$ be respectively the unique solution to 
\begin{equation}\label{eq-uu-delta}
\Div (s_\delta A \nabla u_\delta) + s_0 k^2 \Sigma u_\delta = f \mbox{ in } \Omega, 
\end{equation}
and 
\begin{equation}\label{eq-uu}
\Delta u + k^2 u = f \mbox{ in } \Omega. 
\end{equation}
As in \cite{Ng-Complementary}, we assume that
\begin{equation}
	\mbox{equation}\, \eqref{eq-uu}\, \mbox{with } f = 0, \mbox{has only zero solution in}\, H^1_0(\Omega).
\end{equation}

Our result on cloaking using complementary media is: 

\begin{theorem}\label{thm1} Let  $d=2, \, 3$, $f \in L^{2}(\Omega)$ with $\supp f \subset \Omega \setminus  B_{r_{3}}$ and let  $u$ and  $u_\delta$ in $ H^1_{0}(\Omega)$ be the unique solution to \eqref{eq-uu-delta} and \eqref{eq-uu} resp. There exists $\gamma_0> 1$, depending {\bf only} on $\Lambda$ and the Lipschitz constant of $\hat a$  such that if $1 < \gamma < \gamma_0 $ then   
\begin{equation}\label{key-point}
u_{\delta} \to u \mbox{ weakly in } H^{1} (\Omega \setminus B_{r_3}) \mbox{ as } \delta \to 0. 
\end{equation}
\end{theorem}
For an observer outside $B_{r_3}$, the medium in $B_{r_3}$ looks like the homogeneous one by \eqref{key-point} (and also \eqref{eq-uu}): one has cloaking.

\begin{remark} \fontfamily{m} \selectfont  The case $k=0$ was established in \cite{Ng-Negative-Cloaking}.
\end{remark}

The proof of Theorem~\ref{thm1} is given in Section~\ref{sect-cloaking}.  It is based on the removing localized singularity technique introduced in  \cite{Ng-Negative-Cloaking} and uses a new three sphere inequality (Theorem~\ref{thm-threesphere}) discussed in the next section. The discussion on illusion optics is given in Section~\ref{sect-illusion}.

\section{Three spheres inequalities}\label{Sect-ThreeSphereInequality}
 
 Let $v$ be an holomorphic function defined in  $B_{R_3}$, Hadamard in \cite{Hadamard} proved 
 the following famous three spheres inequality:
 \beq
	\|v\|_{L^{\infty}(\pB_{R_2})} \leq \|v\|_{L^{\infty}(\pB_{R_1})}^{\alpha} \|v\|_{L^{\infty}(\pB_{r_3})}^{1 - \alpha}
\eeq{Hadamard}  
for all $0 < R_1 < R_2 < R_3$ where 
\begin{equation*} 
\dsp 	\alpha = \log \Big( \frac{R_3}{R_2} \Big) \Big/ \log \Big(\frac{R_3}{R_1} \Big).
\end{equation*}
A three spheres inequality for general elliptic equations was proved by Landis \cite{Landis} using Carleman type estimates. 
 Landis proved \cite[Theorem 2.1]{Landis} that \footnote{In fact, \cite[Theorem 2.1]{Landis} deals with the non-divergent form; however since $M$ is assumed $C^2$, the two forms are equivalent.} if $v$ is a solution to 
\begin{equation}\label{Equation}
\dive (M \nabla v) +  \vec{b} \cdot \nabla v + c v  = 0 \mbox{ in } B_{R_3},  
\end{equation}
where $M$ is elliptic, symmetric,  and of class $C^2$,  $\vec b, c \in C^1$, {\bf and  $c \le 0$}, then 
there is a constant $C>0$ such that
\beq
	\|v\|_{L^{\infty}(\pB_{R_2})} \leq C \|v\|_{L^{\infty}(\pB_{R_1})}^{\alpha} \|v\|_{L^{\infty}(\pB_{r_3})}^{1 - \alpha}
\eeq{Landis}  
for some $\alpha  \in (0, 1)$ depending only on $R_2/R_1, R_2/R_3$, the ellipticity constant of $M$,  and the regularity constants of $M$, $b$, and $c$. {\bf The assumption $c \le 0$ is crucial} and this is discussed in the next paragraph.  Another proof was obtained by Agmon \cite{Agmon} in which he used the logarithmic convexity. Garofalo and Lin in \cite{GarofaloLin} established similar results where the $L^\infty$-norm is replaced by the $L^2$-norm, and $M$ is of class $C^1$, $\vec{b}$ and $c$ are in $L^\infty$: 
\beq
	\|v\|_{L^{2}(\pB_{R_2})} \leq C \|v\|_{L^{2}(\pB_{R_1})}^{\alpha} \|v\|_{L^2(\pB_{r_3})}^{1 - \alpha}
\eeq{Lin}  
using  the frequency function. 
  
A typical example of \eqref{Equation} when $c>0$ is the  Helmholtz equation: 
\begin{equation}\label{H}
\Delta v +  k^2 v = 0 \mbox{ in } B_{R_3}.   
\end{equation}
Given $k> 0$, neither \eqref{Lin} nor \eqref{Landis} holds for all $R_1 < R_2 < R_3$. Indeed, first consider the case $d=2$. 
It is clear that  for $n \in \mZ \setminus \{0\}$, the function $J_n(k r ) e^{i n \theta}$  is a solution to \eqref{H}  in $\mR^2 \setminus \{0\}$, where $J_n$ is  the  Bessel function of order $n$. By taking $R_1$, $R_2$, and $R_3$ such that $J_n(kR_1) = 0 \neq J_n(kR_2)$, one reaches the fact that neither \eqref{Lin} nor \eqref{Landis}  is valid. The same conclusion holds in the higher dimensional case by similar arguments. In the case $c > 0$,  \eqref{Lin}  holds  {\bf under the smallness of  $R_3$} (see e.g.,  \cite[Theorem 4.1]{AlessandriniRondi}); this condition is equivalent to the smallness of  $c$ for a fixed $R_3$ by a scaling argument.  

In this paper, we establish  a new type of  three spheres inequalities without imposing  the smallness condition on $R_3$. This inequality will play an important role in the proof of Theorem~\ref{thm1}. 
Define 
\begin{equation}\label{H-norm}
\| v \|_{\bH(\pB_r)} = \| v\|_{H^{1/2}(\pB_r)} + \| M \nabla v \cdot \nu \|_{H^{-1/2}(\pB_r)}. 
\end{equation} 
Here and in what follows, $\nu$ denotes the outward normal vector on a sphere. 

\medskip
Our result on three spheres inequalities is:

\begin{theorem}\label{thm-threesphere}
	Let $d \ge 2$, $c_1, c_2 >0$, $0 < R_* <  R_1  < R_2 <  R_3 < R^*$, and let $M$ be  a Lipschitz  uniformly elliptic symmetric matrix-valued function defined in $B_{R^*}$. 
Assume $v \in H^1(\annulus{R_3}{R_1})$ satisfies
	\begin{equation}\label{eq-Main}
|\Div (M \nabla v)| \le  c_1 |\nabla v|  + c_2|v|, \quad \mbox{ in } \X.
\end{equation}
Then, for any $\lambda_0 > 1$ with $R_2 \in (\lambda_ 0 R_1, R_3/\lambda_0),$ there exist a constant $C$ and $q \ge 1$, depending on the elliptic and the Lipschitz constant of $M$,  $C$ also depends on $c_1,$ $c_2,$  $R_*, R^*,  d$, and $\lambda_0$ such that
 \begin{equation}\label{P1-2}
\| v\|_{\bH(\pB_{R_2})} \le C\| v\|_{\bH(\pB_{R_1})} ^{\alpha} \| v\|_{\bH(\pB_{R_3})}^{1-\alpha}
\end{equation}
where
\begin{equation}\label{def-alpha}
	\alpha := \frac{R_2^{-q} - R_3^{-q}}{R_1^{-q} - R_{3}^{-q}}.
\end{equation}
\end{theorem}

In Theorem ~\ref{thm-threesphere}, one does not impose any smallness condition on $R_1, R_2, R_3$ and the exponent $\alpha$ is independent of $c_1$ and $c_2$.
The proof of Theorem~\ref{thm-threesphere} is  inspired from the approach of Protter in  \cite{Protter60}.  Nevertheless,  different test functions are used.  The ones  in \cite{Protter60} are too concentrated at 0 and not suitable for our purpose.  The connection between  three spheres inequalities and the unique continuation principle, and the application of three spheres inequalities for the stability of Cauchy problems can be found in \cite{AlessandriniRondi}.  

\medskip
The proof of Theorem~\ref{thm-threesphere} is presented in the next two subsections.

%
%
\subsection{Preliminaries}

This section contains several lemmas used in the proof of Theorem~\ref{thm-threesphere}. These lemmas are in the spirit of \cite{Protter60}. Nevertheless, the test functions used here are different from there. Let $0 < R_1 < R_3 < + \infty$.  In this section, we assume  that 
 $M $ is  a Lipschitz symmetric matrix-valued function defined in $\XX$ and satisfies
 \begin{equation*}
\frac{1}{\Lambda} |\xi|^2 \le  M(x) \xi \cdot \xi  \le \Lambda |\xi|^2  \quad \forall \, \xi \in \mR^d,  
\end{equation*}	 
for a.e. $x \in \XX$,  for some $\Lambda \ge 1$.  Set
\begin{equation}\label{def-L1}
L :  = \| M \|_{L^\infty} + R_3\|\nabla  M \|_{L^\infty}. 
\end{equation}

The first lemma is:

\begin{lemma} \label{lem-prepare1} Let $d \ge 2$ and $z \in H^2(\X)$. We have
\begin{equation*}
 \int_\X   (x \cdot M \nabla z) \,   \Div (M\nabla z) \ge   -    \int_{\X} C  L^2|\nabla z|^2   -  \int_{\pX} C L^2 r |\nabla z|^2, 
\end{equation*}
for some positive constant $C$ depending only on $d$. 
\end{lemma}

\noindent{\bf Proof.} An integration by parts gives
\begin{multline}\label{1.4}
	 \int_\X   (x \cdot M \nabla z) \,   \Div (M\nabla z) 
	 = - \int_{\X} \nabla ( x \cdot M  \nabla z) \cdot M\nabla z \\+ \int_{\pX}  (x \cdot M\nabla z) \,  M \nabla z \cdot \nu. 
\end{multline}
Using the symmetry of $M$, we have \footnote{In what follows, the repeated summation is used. }
\begin{equation}\label{1.6}
\frac{\partial }{\partial x_i}(x \cdot M\nabla z) = \frac{\partial }{\partial x_i} \Big( M_{k j} x_j \frac{\partial z}{\partial x_k} \Big) = M_{k j} x_j \frac{\partial^2 z}{\partial x_i \partial x_k} + M_{k i } \frac{\partial z}{\partial x_k}  + x _j \frac{\partial M_{k j}}{\partial x_i}   \frac{\partial z}{\partial x_k}. 
\end{equation}
and
\begin{multline}\label{1.6-1}
		-\int_{\X} 2 x_j M_{kj} \frac{\partial^2 z}{\partial x_i \partial x_k}M_{il} \frac{\partial z}{\partial x_l}
		 = - \int_\X x_j M_{kj} M_{il} \frac{\partial}{\partial x_k}\left(\frac{\partial z}{\partial x_i} \frac{\partial z}{\partial x_l}\right) \\
		= \int_\X \frac{\partial (x_j M_{kj} M_{il})}{\partial x_k}\frac{\partial z}{\partial x_i} \frac{\partial z}{\partial x_l} - \int_{\pX}  x_j M_{kj} M_{il}\frac{\partial z}{\partial x_i} \frac{\partial z}{\partial x_l} \nu_k. 
\end{multline}
We derive from \eqref{1.6} and \eqref{1.6-1} that
\begin{equation}
- \int_{\X} \nabla ( x \cdot M\nabla z) \cdot M\nabla z \\
 \ge    -    \int_{\X} C  L^2|\nabla z|^2  -   \int_{\pX}  C L^2 r|\nabla z|^2.  
\label{1.7}
\end{equation}
The conclusion now follows from  \eqref{1.4} and  \eqref{1.7}.  \proofend

\medskip 
The second lemma is 

\begin{lemma}\label{lem-prepare2} Let $d \ge 2$, $\beta \in \mR$, and $z \in H^2(\X)$. There exists  $p_{\Lambda, L} \geq 1$ such that if $p \ge p_{\Lambda, L}$ and $|\beta| R_3^{-p} \ge 2$ then
\begin{multline*}
 \int_\X    e^{\beta r^{-p}} (M x \cdot \nabla |z|^2) \,  \Div (M\nabla  e^{-\beta r^{-p}}) +  \int_{\pX} C L^2 p^2 \beta^2  r^{- 2 p-1} |z|^2 \\[6pt]
 \ge  \int_\X  \frac{1}{2}\Lambda^{-2} p^3 \beta^2 r^{-2p - 2}|z|^2 - \int_\X C L^2 |\nabla z|^2. 
\end{multline*}
for some positive constant $C$ depending only on $d$. 
\end{lemma}

\noindent{\bf Proof.} A computation yields 
\begin{align*}
	\Div (M\nabla \testfn) =   p \beta  \testfn \big[p \beta r^{-2 p - 4} - ( p + 2) r^{-p - 4}\big] x \cdot Mx +p \beta r^{-p - 2} \testfn \, \Div (Mx).
		\end{align*}
An integration by parts gives
\begin{equation}\label{1.9}
	 \int_\X   e^{\beta r^{-p}} (M x \cdot \nabla |z|^2) \,   \Div (M\nabla \testfn) =  P + Q. 
\end{equation}	
Here 
\beq
	P = P_1 + P_2 + P_3 \nonumber
\eeq{} with
\begin{equation*}\left\{
\begin{array}{rl}\dsp P_1 &= \ds- \int_{\X} p^2 \beta^2 |z|^2 \Div \big[ r^{- 2 p - 4} (x \cdot M x) M x \big], \\[6pt]
\dsp P_2 &= \ds \int_{\X} p \beta (p+2) |z|^2 \Div \big[ r^{- p - 4} (x \cdot M x) M x \big], \\[6pt]
\dsp P_3 &=   \ds\int_{\X} 2 p \beta r^{-p -2} \Div(M x) z \nabla z \cdot Mx  \big], 
\end{array}\right.
\end{equation*}
and
\begin{equation*}
Q =  \int_{\pX}  p \beta |z|^2   \Big( \big[p \beta r^{- 2p - 4} -( p + 2) r^{-p - 4} \big] x \cdot M x \Big) M x \cdot \nu. 
\end{equation*}
We next estimate $P$ and $Q$. A computation yields 
\begin{equation*}
- \Div \big[ r^{- 2 p - 4} (x \cdot M x) M x \big] = (2p + 4) (x \cdot Mx)^2  r^{-2p - 6}  -  r^{-2p - 4}   \Div \big[ (x \cdot Mx) Mx\big ]. 
\end{equation*}
This implies 
\begin{equation}\label{R1}
P_1 \ge  \int_\X p^2 \beta^2  r^{-2p - 2}|z|^2 \big[ (2p + 4)\Lambda^{-2}   -  C L^2 \big] . 
\end{equation}
Similarly, 
\begin{equation}\label{R2}
P_2 \ge - \int_\X (p + 2) p |\beta|  r^{-p - 2}|z|^2 \big[ (p + 4)\Lambda^{-2}   +  C L^2 \big]. 
\end{equation}
A combination of \eqref{R1} and \eqref{R2} yields
\begin{equation}\label{R12}
P_1 + P_2 \ge \int_\X  \Lambda^{-2} p^3 \beta^2 r^{-2p - 2}|z|^2. 
\end{equation}
Here we used the fact that $p \ge p_{\Lambda, L}$ and $|\beta| R_3^{-p} \ge 2$. 
On the other hand, using Cauchy's inequality, we have
\begin{align*}
		|P_3| \leq \int_{\X} p^2 \beta^2 r^{-2p - 2} L^2 |z|^2 + \int_{\X} C L^2 |\nabla z|^2.
\label{R}
\end{align*}
It follows from  \eqref{R12} that 
\begin{align}
	P &\geq \int_\X \frac{1}{2} p^3\beta^2 \Lambda^{-2} r^{-2p - 2}|z|^2 - \int_{\X}  C L^2 |\nabla z|^2,
\end{align}
provided that $p \geq 2\Lambda^{2} L^2.$ Since 
\begin{equation*}
|Q| \le   \int_{\pX}  2 \Lambda^2 p^2 \beta^2 r^{-2p -1 } |z|^2.  
\end{equation*}
the conclusion follows. \proofend

\medskip 
Using Lemmas~\ref{lem-prepare1} and \ref{lem-prepare2}, we can prove the following result.

\begin{lemma}\label{lem1} Let $d \ge 2$, $\beta \in \mR$,  and $v \in H^2(\X)$. There exists a positive constant $p_{\Lambda, L} \geq 1$ such that if $p \ge p_{\Lambda, L}$ and $|\beta| R_3^{-p} \ge 2$ then
\begin{multline*}
 \int_\X \frac{r^{p + 2}e^{2\beta r^{-p}}}{2p |\beta|}\big[\Div (M\nabla v)\big]^2 +  
  \int_{\X}  C  L^2  e^{2 \beta r^{-p}}  |\nabla v|^2 \\
  + \int_{\pX} CL^2 p^2 \beta^2  r^{-2p-1} e^{2 \beta r^{-p}}|v|^2  +  \int_{\pX} C L^2 r e^{2 \beta r^{-p}}|\nabla v|^2 \\
\ge \int_\X \frac{1}{2} \Lambda^{-2} p^3 \beta^2 r^{-2p - 2} e^{2 \beta r^{-p}} |v|^2, 
\end{multline*}
for some positive constant $C$ depending only on $d$. 
\end{lemma}


\noindent {\bf Proof.} By considering the real part and the imaginary part of $v$ separately, one might assume that $v$ is real. 
Set
\begin{equation*}
	z = e^{\beta r^{-p}} v \quad \mbox{ equivalently }  v = e^{-\beta r^{-p}} z.
\end{equation*}
Since  $\Div \big(M \nabla (g h) \big) = 2 \nabla h \cdot  M\nabla g + h \Div(M \nabla g) + g \Div (M \nabla h)$ ($M$ is symmetric), it follows that 
\begin{align*}
	\Div (M\nabla v) = 2\beta p r^{-p  - 2}  e^{-\beta r^{-p}} x \cdot M\nabla z  +  e^{-\beta r^{-p}} \Div (M\nabla z)  + z \Div (M\nabla  e^{-\beta r^{-p}}).
\end{align*}
Using the inequality $(a + b + c)^2 \geq 2 a (b + c)$,  we obtain
\begin{equation*}
\frac{1}{2 }\big[ \Div (M\nabla v)  \big]^2 \ge 2|\beta| p r^{-p  - 2}  e^{-\beta r^{-p}} (x \cdot M\nabla z) \Big( e^{-\beta r^{-p}} \Div (M\nabla z)  + z \Div( M\nabla  e^{-\beta r^{-p}})\Big). 
\end{equation*}
This implies 
\begin{multline*}
\int_\X \frac{r^{p + 2}e^{2\beta r^{-p}}}{2p|\beta|}\big[\Div (M\nabla v)\big]^2 
\geq \int_\X 2   (x \cdot M\nabla z ) \,  \Div (M\nabla z) \\+ \int_\X    e^{\beta r^{-p}} (Mx \cdot \nabla |z|^2) \,  \Div (M\nabla  e^{-\beta r^{-p}}). 
\end{multline*}
Applying Lemmas~\ref{lem-prepare1} and \ref{lem-prepare2},  we have 
\begin{multline}\label{p1-Lem3}
 \int_\X \frac{r^{p + 2}e^{2\beta r^{-p}}}{2p|\beta|}\big[\Div (M\nabla v)\big]^2 
\ge \int_\X \Big( \Lambda^{-2} p^3 \beta^2 r^{-2p - 2}|z|^2 
-    C L^2 |\nabla z|^2 \Big) \\[6pt] - \int_{\pX}  \Big( C L^2 p^2 \beta^2  r^{-2 p-1} |z|^2 + C L^2 r |\nabla z|^2 \Big).
\end{multline}
Since $z = e^{\beta r^{-p}} v $, 
\begin{equation}\label{p2-Lem3}
 |\nabla z|^2 \leq 2 e^{2 \beta r^{-p}} (|\nabla v|^2 + p^2\beta^2 r^{-2 p - 2}|v|^2). 
\end{equation}
A combination of \eqref{p1-Lem3} and \eqref{p2-Lem3} yields, since $p \ge p_{\Lambda, L}$, 
\begin{multline*}
 \int_\X \frac{r^{p + 2}e^{2\beta r^{-p}}}{2p|\beta|}\big[\Div (M\nabla v)\big]^2 
\ge \int_\X e^{2 \beta r^{-p}} \left( \frac{1}{2} \Lambda^{-2} p^3 \beta^2 r^{-2p - 2}  |v|^2  - 
C L^2    |\nabla v|^2 \right) \\[6pt]- \int_{\pX} e^{2 \beta r^{-p}} \left( C L^2 p^2 \beta^2  r^{-2p-1} |v|^2   +  C L^2 r |\nabla v|^2  \right).
\end{multline*}
The conclusion follows. \proofend

\medskip 
We also have

\begin{lemma}\label{lem2}  Let $d \ge 2$, $\beta \in \mR$, and $v \in H^2(\X)$. 
There exists a positive constant $p_{\Lambda, L} \geq 1$ such that if $p \ge p_{\Lambda, L}$ and $|\beta| R_3^{-p} \ge 2$ then
\begin{multline*}
 \int_\X e^{2\beta r^{-p}} v \,  \Div (M\nabla v) +  \int_\X   e^{2 \beta r^{-p}} |\nabla v|^2 \\
 \le \int_\X  C \beta^2 p^2 r^{-2p - 2} e^{2 \beta r^{-p}}  |v|^2+ \int_{\pX} C e^{2\beta r^{-p}} (r|\nabla v|^2 +  r^{-1}|v|^2), 
\end{multline*}
for some positive constant $C$ depending only on $d$, $\Lambda$, and $L$. 
\end{lemma}

\noindent{\bf Proof.}
We have
\begin{equation}\label{s1}
- \int_\X e^{2\beta r^{-p}} v \,  \Div (M\nabla v)  = \int_\X M\nabla v \cdot \nabla (e^{2\beta r^{-p}} v) - \int_{\pX} e^{2\beta r^{-p}} v M \nabla v \cdot \nu. 
\end{equation} 
On the other hand,  
\begin{equation}\label{s2}
 \int_\X M\nabla v \cdot \nabla (e^{2\beta r^{-p}} v) 
	=\int_\X  \Big(e^{2\beta r^{-p}} M\nabla v \cdot \nabla v - 2 \beta p  r^{-p - 2}e^{2\beta r^{-p}} v M\nabla v  \cdot x \Big)
\end{equation} 
and 
\begin{equation}\label{s3}
 \int_{\pX} e^{2\beta r^{-p}} v M \nabla v \cdot \nu \le  \int_{\pX} e^{2\beta r^{-p}} (r|\nabla v|^2 + L^2 r^{-1}|v|^2). 
\end{equation}
Since
\begin{equation*}
2  \beta p  r^{-p - 2}v M\nabla v \cdot x \le \frac{1}{2} \Lambda^{-1} |\nabla v|^2  + 8 \beta^2 p^2 L^2 \Lambda r^{-2p - 2} |v|^2, 
\end{equation*}
we derive from  \eqref{s1}, \eqref{s2}, and \eqref{s3} that 
\begin{multline*}
 \int_\X e^{2\beta r^{-p}} v \,  \Div (M\nabla v) +  \int_\X \frac{1}{2} \Lambda^{-1} e^{2 \beta r^{-p}} |\nabla v|^2 \\
 \le \int_\X  C \beta^2 p^2 r^{-2p - 2} e^{2 \beta r^{-p}}  |v|^2+ \int_{\pX} C e^{2\beta r^{-p}} \big(r|\nabla v|^2 +  r^{-1}|v|^2 \big). 
\end{multline*}
The conclusion follows. \proofend

\medskip 

Combining the inequalities of Lemmas \ref{lem1} and \ref{lem2}, we obtain
\begin{lemma}\label{lem-main}
	 Let $d \ge 2$, $\beta \in \mR$, and  $v \in H^2(\X)$. There exists a positive constant $p_{\Lambda, L} \geq 1$ such that if   $p \ge p_{\Lambda, L}$ and  $|\beta| \geq 2 R_3^{-p}$ then
\begin{multline}\label{comb}
	\int_{\X} \testfnn  |\beta| p \Big(p^3 \beta^2 r^{-2p - 2} e^{2 \beta r^{-p}} |v|^2 +  |\nabla v|^2 \Big) \leq C \int_\X r^{p+2} \testfnn|\Div(M\nabla v)|^2   \\+ C \int_{\pX}   |\beta| p e^{2\beta r^{-p}}  \big(r|\nabla v|^2 + p^2 \beta^2  r^{-2p -1}|v|^2 \big).
\end{multline}
for some positive constant $C$ depending only on $d$, $\Lambda$, and $L$. 
\end{lemma}

\noindent{\bf Proof.}  Note that 
\begin{align*}
|v \dive (M \nabla v)| \le  p|\beta|^2 |v|^2 r^{- 2p -2 } + \frac{4 }{p |\beta|} |\dive (M \nabla v)|^2 r^{p+2}. 
\end{align*}
The conclusion now follows from Lemmas \ref{lem1} and \ref{lem2}.  The details are left to the reader. 
\proofend 

\subsection{Proof of Theorem~\ref{thm-threesphere}}

Let
\begin{equation*}
1 < \lambda < \lambda_0, 
\end{equation*}
(which will be defined later) and set 
\begin{equation*}
D = \annulus{\lambda R_3}{R_1/\lambda}. 
 \end{equation*}
Let $u_1 \in H^1(D \setminus \partial B_{R_1})$ and $u_3 \in H^1(D \setminus \partial B_{R_3})$ be respectively the unique solution to 
\begin{equation*} \left\{
\begin{array}{cl}
	\Div(M\nabla u_1) = 0 &\mbox{in } D \setminus \partial B_{R_1},\\[6pt]
	\ds [u_1] = v; \; \ds [M\nabla u_1 \cdot \nu] = M\nabla v \cdot \nu &\mbox{on } \pB_{R_1},\\[6pt]
	u_1 = 0 &\mbox{on } \partial D, 
\end{array} \right.
\end{equation*}
and 
\begin{equation*} \left\{
\begin{array}{cl}
	\Div(M\nabla u_3) = 0 &\mbox{in } D \setminus B_{R_3},\\[6pt]
	\ds [u_3] = v; \; \ds [M\nabla u_3 \cdot \nu] = M\nabla v \cdot \nu &\mbox{on } \pB_{R_3},\\[6pt]
	u_3 = 0 &\mbox{on } \partial D.
\end{array} \right.
\end{equation*}
Here and in what follows, $[\cdot]$ denotes the jump across a sphere and  $\nu$ denotes the unit outward normal vector on a sphere. 
It follows that 
\begin{equation}\label{inf-R1-Thm1}
\| u_1 \|_{H^1(D\setminus \partial B_{R_1} )} \le C \|v \|_{\bH(\partial B_{R_1})}, \quad 
	\|u_1\|_{H^{3/2}(\pB_{R_1/\gamma})} \leq C\|v\|_{\bH(\pB_{R_1})}
\end{equation}
and 
\begin{equation}\label{inf-R3-Thm1}
\| u_3 \|_{H^1(D \setminus \partial B_{R_3})} \le C \|v \|_{\bH(\partial B_{R_3})}, \quad 	\|u_3\|_{H^{3/2}(\pB_{\gamma R_3})} \leq C\|v\|_{\bH(\pB_{R_3})}.
\end{equation}
Here and in what follows in this proof, $C$ denotes a positive constant depending only on the elliptic and the Lipschitz constant of $M$, $c_1$, $c_2$,  $\lambda_0$, $R_*$, $R_*$, and $d$. 
Set 
\begin{equation*}
d_1 = (\lambda -1)R_1 \quad \mbox{ and } \quad d_3 = (\lambda - 1)R_3/\lambda. 
\end{equation*}
Let $\varphi_1, \; \varphi_3 \in C^2_{c}(\mR^d)$ be such that 
\begin{equation*}
	\varphi_1 = \left\{
		\begin{array}{ll}
			1 & \mbox{ in } B_{R_1 + d_1/3} \setminus B_{R_1}, \\[6pt]
			0 &\mbox{in } \mR^d \setminus (B_{R_1 + d_1/2} \setminus B_{R_1/\lambda})
		\end{array}
	\right. \quad \mbox{and } \quad 
	\varphi_3 = \left\{
		\begin{array}{ll}
			1 &\mbox{in } B_{R_3} \setminus  B_{R_3-d_3/3},\\[6pt]
			0 &\mbox{in } \mR^d \setminus (B_{\lambda R_3} \setminus B_{R_3 - d_3/2}).
		\end{array}
	\right.
\end{equation*}
Define 
\begin{equation}\label{def-V-Thm1}
V = 
\left\{\begin{array}{cl}
 v- \varphi_1 u_1 - \varphi_3 u_3  & \mbox{ in } \X,  \\[6pt]
 - \varphi_1 u_1 - \varphi_3 u_3  & \mbox{ in } D \setminus (\X). 
 \end{array}\right.
\end{equation}
Applying Lemma~\ref{lem-main}, we obtain, for $|\beta| > 2(\gamma R_{3})^{p}$,  
\begin{equation}\label{2.24}
	C \int_{D} \testfnn \beta (\beta^2  |V|^2 +    |\nabla V|^2) \leq  \int_{D}  \testfnn |\Div(M\nabla V)|^2   + \int_{\partial D}   |\beta|  e^{2\beta r^{-p}}  (|\nabla V|^2 + \beta^2 |V|^2).
\end{equation}

The proof is now quite standard and  divided into two cases.

\medskip 
\noindent \underline{Case 1}: $\| v\|_{\bH(\pB_{R_1})} \le \| v\|_{\bH(\pB_{R_3})}$. 
We deduce from \eqref{2.24} that for $\beta \ge \beta_0$, 
\begin{equation}
	C \int_{D} \testfnn (|V|^2 +  |\nabla V|^2) \leq     \beta^2  e^{2\beta  \hat R_{3}^{-p}} \| v\|_{\bH(\partial B_{R_3})}^2 
	+  \beta^2 e^{2\beta \hat R_{1}^{-p}} \|v\|_{\bH(\partial B_{R_1})}^2, 
\label{2.23}
\end{equation}
where 
\begin{equation*}
\hat R_3 = R_3/ \lambda \quad \mbox{ and } \quad \hat R_1 = R_1 / \lambda. 
\end{equation*}
This implies
\begin{equation}\label{point-point}
C \|v \|_{\bH(\partial B_{R_2})} \leq     \beta  e^{\beta (\hat R_3^{-p}  - R_2^{-p})} \| v\|_{\bH(\partial B_{R_3})}
	+  \beta  e^{\beta (\hat R_1^{-p} - R_2^{-p})} \|v\|_{\bH(\partial B_{R_1})}. 
\end{equation} 
Define  $\alpha' \in (0, 1)$ and $\beta> 0$ as follows
\begin{equation*}
\alpha' = \frac{R_2^{-p} - \hat R_3^{-p}}{\hat R_1^{-p} - \hat R_3^{-p}}
\quad \mbox{ and }
\beta (R_2^{-p} - \hat R_1^{-p}) =  (1 - \alpha') \ln \Big(\|v \|_{\bH(\pB_{R_3})} \big/\| v\|_{\bH(\pB_{R_1})} \Big) \;  \footnote{Here we assume that $\| v\|_{\bH(\pB_{R_1})} \neq 0$ since otherwise $v= 0$. This fact is a consequence of the unique continuation principle and can be obtained from \eqref{point-point} by letting $\beta \to \infty$.}.
\end{equation*}
 Note that $0 < \alpha' < 1 $ since $R_2 < R_3/ \gamma$. 
We assume that $\|v \|_{\bH(\pB_{R_3})}> C \| v\|_{\bH(\pB_{R_1})}$ for some large $C$ such that $\beta \ge \max\{2R_3^{-p}, 2, \beta_0\}$ since if $\|v \|_{\bH(\pB_{R_3})} <  C \| v\|_{\bH(\pB_{R_1})}$, the conclusion holds for any $\alpha \in (0, 1)$ by taking $\beta = \max\{2 R_3^{-p}, 2, \beta_0\}$ in \eqref{point-point}. 
It follows from \eqref{point-point} and the choice of $\alpha'$ and $\beta$ that 
\begin{equation}\label{with beta}
\| v \|_{\bH(\pB_{R_2})}
\le C \beta \|v \|_{\bH(\pB_{R_1})}^{\alpha'}\| v\|_{\bH(\pB_{R_3})}^{1-\alpha'}. 
\end{equation}
Define
\beq
	\alpha: = \frac{R_2^{-2p} - R_3^{-2p}}{R_1^{-2p} - R_3^{-2p}}. 
\eeq{alpha}
It is clear that $\alpha < \frac{R_2^{-p} - R_3^{-p}}{R_1^{-p} - R_3^{-p}}$. Hence, by choosing $\lambda$ close to 1,  
\begin{equation}\label{LL}
\alpha < \alpha' \quad \forall R_2 \in (\gamma_0 R_1, R_3/ \gamma_0). 
\end{equation}
A combination of  \eqref{with beta} and \eqref{LL} implies
\begin{equation}\label{2.25}
\| v \|_{\bH(\pB_{R_2})}
\le C \|v \|_{\bH(\pB_{R_1})}^{\alpha}\| v\|_{H(\pB_{R_3})}^{1-\alpha}.
\end{equation}

\medskip 
\noindent \underline{Case 2}: $\| v\|_{\bH(\pB_{R_1})} \ge \| v\|_{\bH(\pB_{R_3})}$. The proof is similar to the previous case by  considering  $\beta < - 2 (\gamma R_{3})^{-p}$. The details are left to the reader. The proof is complete. \proofend

\section{Cloaking using complementary media. Proof of Theorem~\ref{thm1}} \label{sect-cloaking}

This section containing three subsections is devoted to the proof of Theorem~\ref{thm1}. In the first subsection, we present two useful lemmas. The proof of Theorem~\ref{thm1} is given in the second subsection. 

\subsection{Preliminaries}

In this section, we present two lemmas which will be used in the proof of Theorems~\ref{thm1} and \ref{thm2}.  The first lemma is on a change of variables and follows from  \cite[Lemma 1]{Ng-Complementary}. 
\begin{lemma}\label{lem-TO} Let $d \ge 2$, $k>0$,  $0 < R_1 < R_2 < R_3$ with $R_3 = R_2^2/ R_1$, $a \in [L^\infty(B_{R_3} \setminus \overline B_{R_2})]^{d \times d}$ be a matrix valued function, $\sigma \in L^\infty(B_{R_3} \setminus \overline B_{R_2})$ be a complex function, 
and $K:B_{R_2} \setminus \bar B_{R_1} \to B_{R_3} \setminus \bar B_{R_2}$ be the Kelvin transform with respect to $\partial B_{R_2}$, i.e., 
\begin{equation*}
K(x) = R_2^2 x/ |x|^2. 
\end{equation*}
For $v \in H^1(\annulus{R_3}{R_2})$, define $w = v \circ K^{-1}$. Then 
\begin{equation*}
\dive (a \nabla v) + k^2 \sigma v = 0 \mbox{ in }  \annulus{R_3}{R_2}
\end{equation*}
if and only if
\begin{equation*}
\dive (K_*a \nabla w)  + k^2 K_* \sigma w = 0 \mbox{ in } \annulus{R_2}{R_1}. 
\end{equation*}
Moreover, 
\begin{equation*}
w = v \quad \mbox{ and } \quad K_*a \nabla w \cdot \nu = - a \nabla v \cdot \nu  \mbox{ on } \partial B_{R_2}.    
\end{equation*}
\end{lemma}

The second lemma is a stability estimate for solutions of \eqref{eq-uu-delta}. 

\begin{lemma}\label{lem-stability1} Let $0 < \delta < 1$, $f \in L^2(\Omega)$, and let  $A \in [L^\infty(\Omega)]^{d \times d}$ and $\Sigma \in L^\infty(\Omega, \mathbb{C})$ be such that $A$ is Lipschitz and uniformly elliptic and $\Im(\Sigma) \ge 0$. Assume that 
$u_\delta \in H^1_0(\Omega)$ is  the unique solution to \eqref{eq-uu-delta}. Then
\begin{equation}\label{est-stability}
\| u_\delta\|_{H^1(\Omega)}^2 \le C \left(\delta^{-1} \|f\|_{L^2(\Omega)}\|u_{\delta}\|_{L^2({\textrm supp} f)} + \|f\|_{L^2(\Omega)}^2\right),  
\end{equation}
for some positive constant $C$ independent of $\delta$ and $f$. 
\end{lemma}

Lemma~\ref{lem-stability1} is a variant of  \cite[Lemma 1]{Ng-Complementary}. The case $k=0$ and its variant in the case $k>0$ were considered in  \cite{Ng-Negative-Cloaking} and \cite{Ng-Superlensing} respectively. The proof is similar to the one of \cite[Lemma 1]{Ng-Complementary}. For the convenience of the reader, we present the proof. 

\medskip
\noindent{\bf Proof.}  The existence  and uniqueness of $u_\delta$ are given in \cite{Ng-Complementary}.  We only establish \eqref{est-stability} by contradiction. Assume that \eqref{est-stability} is not true. Then 
there exist $\delta_n \to 0$ and  $(f_n) \subset L^2(\Omega)$ such that  
\begin{equation}\label{contradict-assumption}
\| u_n\|_{H^1(\Omega)} =1 \mbox{ and } \frac{1}{\delta_n} \|f_n\|_{L^2(\Omega)} \| u_n\|_{L^2({\textrm supp} f_n)} + \| f_n\|_{L^2(\Omega)}^2  \to 0,
\end{equation}
as $n \to \infty$, where $u_n \in H^1_0(\Omega)$ is the unique solution to 
\begin{equation}\label{eq-v-delta}
\dive (\epss_{\delta_n} A \nabla u_n) + k^2 \epss_0 \Sigma u_n= f_n  \mbox{ in } \Omega.
\end{equation}
Without loss of generality, one may assume that $u_n \to u$ weakly in $H^1(\Omega)$ and strongly in $L^2(\Omega)$; moreover, $u \in H^1_0(\Omega)$ and $u$ satisfies
\begin{equation}\label{eq-uu-1}
\dive(s_0 A \nabla  u) + k^2 s_0 \Sigma u = 0 \mbox{ in } \Omega. 
\end{equation}
Multiplying equation \eqref{eq-v-delta} by $\bar u_n$ (the conjugate of $u_n$) and integrating on $\Omega$, we have
\begin{equation*}
\int_{\Omega} \epss_{\delta_n}  A \nabla u_n \cdot \nabla \bar u_n   \, dx -  \int_{\Omega} k^2 \epss_0 \Sigma |u_n|^2 \, dx = - \int_\Omega f_n \bar u_n \, dx.
\end{equation*}
Considering the imaginary part and using the fact that
\begin{equation*}
\frac{1}{\delta_n}\Big| \int_\Omega f_n \bar u_n \Big| \le \frac{1}{\delta}\|f_n \|_{L^2(\Omega)} \| u_n\|_{L^2({\textrm supp} f_n)} \to 0 \mbox{ as $n \to \infty$ by }\eqref{contradict-assumption},
\end{equation*}
we obtain, by \eqref{cond-a},
\begin{equation}\label{toto1}
\| \nabla u_n\|_{L^2(\annulus{r_2}{r_1})} \to 0 \mbox{ as } n \to \infty. 
\end{equation}
Since $\dive(A \nabla u_n) + k^{2} \Sigma u_n = f_n$ in $B_{r_2} \setminus B_{r_1}$ and $f_n \to 0$ in $L^2(\Omega)$, it follows from  that $u_n \to 0$ in the distributional sense. This in turn implies  
\begin{equation}\label{toto1-1}
\| u_n \|_{L^{2}(B_{r_2} \setminus B_{r_1})} \to 0 \mbox{ as } n \to \infty.   
\end{equation}
A combination of \eqref{toto1} and \eqref{toto1-1} yields
\begin{equation}\label{part-1}
\| u_n \|_{H^{1}(B_{r_2} \setminus B_{r_1})} \to 0  \mbox{ as } n \to \infty. 
\end{equation}
Hence   
\begin{equation*}
u= 0 \mbox{ in } B_{r_2} \setminus B_{r_1},  
\end{equation*}
and
\begin{equation}\label{limit-1}
 \| u_n \|_{H^{1/2}(\partial B_{r_2})}  + \| u_n \|_{H^{1/2}(\partial B_{r_1})} + \| A \nabla u_n \cdot \nu \|_{H^{-1/2}(\partial B_{r_2})}  + \| A \nabla u_n \cdot \nu\|_{H^{-1/2}(\partial B_{r_1})} \to 0  \mbox{ as }n \to \infty.  
\end{equation}
Since $u=0$ in $B_{r_2} \setminus B_{r_1}$ and $u$ satisfies \eqref{eq-uu-1}, it follows from the unique continuation principle that $u = 0$ in $\Omega$. Hence, since $u_n \to u$ in $L^2(\Omega)$,
\begin{equation}\label{u0}
u_n \to 0 \mbox{ in } L^2(\Omega) \mbox{ as } n \to \infty.
\end{equation}
Multiplying \eqref{eq-v-delta} by $\bar u_n$ and integrating on $\Omega \setminus B_{r_2}$, we have 
\begin{equation*}
\int_{\Omega \setminus B_{r_2}}  A \nabla u_n \cdot \nabla \bar u_n   \, dx -  \int_{\Omega \setminus B_{r_2}} k^2 \epss_0 \Sigma |u_n|^2 \, dx = - \int_\Omega f_n \bar u_n \, dx + \int_{\partial B_{r_2}} A \nabla u_n \cdot \nu \; \bar u_n.
\end{equation*}
Using \eqref{limit-1} and \eqref{u0}, we obtain 
\begin{equation}\label{part-2}
\|\nabla u_n \|_{L^2(\Omega \setminus B_{r_2})} \to 0 \mbox{ as } n \to \infty. 
\end{equation}
Similarly, 
\begin{equation}\label{part-3}
\|\nabla u_n \|_{L^2(B_{r_1})} \to 0 \mbox{ as } n \to \infty. 
\end{equation}
A combination of \eqref{part-1}, \eqref{u0}, \eqref{part-2}, and \eqref{part-3}, we obtain 
\begin{equation*}
\|u_n \|_{H^1(\Omega)} \to 0 \mbox{ as } n \to \infty. 
\end{equation*}
which contradicts  \eqref{contradict-assumption}.  The proof is complete. \proofend


\subsection{Proof of Theorem~\ref{thm1}} 
We use the approach in \cite{Ng-Negative-Cloaking} with some modifications from \cite{Ng-Superlensing} so that the same proof also give the result on illusion optics (Theorem~\ref{thm2} in Section~\ref{sect-illusion}). However, instead of applying the standard three sphere inequality as in \cite{Ng-Negative-Cloaking}, we use Theorem~\ref{thm-threesphere}. 

 We have, by Lemma~\ref{lem-stability1}, 
\begin{equation}\label{haha0}
\| u_{\delta} \|_{H^{1}(\Omega)}^2 \le C \Big( \delta^{-1} \| f\|_{L^2(\Omega)} \|u_{\delta} \|_{L^2(\Omega \setminus B_{r_3})} + \| f\|_{L^2(\Omega)}^2 \Big). 
\end{equation}
As in \cite{Ng-Negative-Cloaking}, let $u_{1, \delta}$ be the reflection of $u_{\delta}$ through $\partial B_{r_{2}}$ by $F$,  i.e., 
\begin{equation}\label{def-u1}
u_{1, \delta } = u_{\delta} \circ F^{-1} \mbox{ in } \mR^{d} \setminus \bar B_{r_{2}}
\end{equation}
and let $u_{2, \delta}$ be the reflection of $u_{1, \delta}$ through $\partial B_{r_2}$ by $G$, i.e., 
\begin{equation}\label{def-u2}
u_{2, \delta } = u_{1, \delta} \circ G^{-1} \mbox{ in }  B_{r_{3}}. 
\end{equation}
By Lemma~\ref{lem-TO}, 
\begin{equation}\label{eqn-u1}
\dive(A \nabla u_{1, \delta}) + \frac{1}{1- i \delta} k^2 \Sigma u_{1, \delta} = 0 \mbox{ in } B_{r_3} \setminus B_{r_2}, 
\end{equation}
\begin{equation}\label{eqn-u2}
\Delta u_{2, \delta} + k^2 u_{2, \delta} = 0 \mbox{ in } B_{r_3}. 
\end{equation}
Applying Lemma~\ref{lem-TO} again and using the fact that $F_*A = A$ in $B_{r_3} \setminus B_{r_2}$, we have  
\begin{equation}\label{TO-1}
u_{1, \delta} = u_{\delta} \Big|_{+}   \mbox{ on } \partial B_{r_{2}} \quad \mbox{ and } \quad  (1 - i \delta ) A \nabla u_{1, \delta} \cdot \nu = A \nabla u_{ \delta} \cdot \nu \Big|_{+} \mbox{ on } \partial B_{r_{2}}. 
\end{equation}
Let $V_{1, \delta} \in H^1(B_{r_3} \setminus B_{r_2})$ be the unique solution to 
\begin{equation} \label{def-V1}
\left\{\begin{array}{cl}
\Div (A \nabla V_{1, \delta}) + k^2 \Sigma V_{1, \delta} = \dsp  - \frac{i \delta }{1 - i \delta} k^2 \Sigma u_{1, \delta}  \mbox{ in } &  B_{r_3} \setminus B_{r_2},  \\[6pt]
A \nabla V_{1, \delta} \cdot \nu- i k V_{1, \delta} = 0 &  \mbox{ on } \partial B_{r_2}, \\[6pt]
V_{1, \delta} = 0 & \mbox{ on } \partial B_{r_3}. 
\end{array}\right. 
\end{equation}
 By Fredholm's theory, 
\begin{equation}\label{V1}
\| V_{1, \delta}\|_{H^1(B_{r_3} \setminus B_{r_2})} \le C \delta \| u_\delta \|_{H^1(\Omega)}. 
\end{equation}
Define  $U_{1, \delta}$ in $B_{r_{3}} \setminus B_{r_2}$ as follows 
\begin{equation}\label{def-bU}
U_{1, \delta} =u_{\delta} - u_{1, \delta} - V_{1, \delta}. 
\end{equation}
Then $U_{1, \delta} \in H^1(B_{r_3} \setminus B_{r_2})$ and $U_{1, \delta}$ satisfies 
\begin{equation*}
\Div (A \nabla U_{1, \delta}) + k^2 \Sigma U_{1, \delta} = 0 \mbox{ in } B_{r_3} \setminus B_{r_2}
\end{equation*}
\begin{equation*}
\| U_{1, \delta}\|_{H^{1/2}(\partial B_{r_2})} + \| A \nabla U_{1, \delta} \cdot \nu \|_{H^{1/2}(\partial B_{r_2})} \le C \delta \| u_\delta \|_{H^1(\Omega)}, 
\end{equation*}
and 
\begin{equation*}
\| U_{1, \delta}\|_{H^{1/2}(\partial B_{r_3})} + \| A \nabla U_{1, \delta} \cdot \nu \|_{H^{1/2}(\partial B_{r_3})} \le C \| u_\delta \|_{H^1(\Omega)}. 
\end{equation*}
Applying  Theorem~\ref{thm-threesphere}, we have  
\begin{equation}\label{U1}
\| U_{1, \delta}\|_{H^{1/2}(\partial B_{\gamma r_2})} + \| A \nabla U_{1, \delta} \cdot \nu \|_{H^{1/2}(\partial B_{\gamma r_2})}  \le C \delta^{\alpha} \| u_\delta\|_{H^1(\Omega)}
\end{equation}
where $\alpha$ is given in \eqref{def-alpha} with $R_1 = r_2$, $R_2 = \gamma r_2$, $R_3 = r_3$. 
By first taking $q$ large enough and then choosing $\gamma_0$ close to 1 enough, from  \eqref{def-alpha}, we can assume that 
\begin{equation}\label{alpha-Thm1}
\alpha > 1/2. 
\end{equation}
Here  is the place where the condition $\gamma < \gamma_0$ is required.  A combination of \eqref{V1} and \eqref{U1} yields 
\begin{equation}\label{part1}
\| u_\delta - u_{1, \delta} \|_{H^{1/2}(\partial B_{\gamma r_2})} + \|A \nabla (u_\delta - u_{1, \delta} ) \cdot \nu\|_{H^{-1/2}(\partial B_{\gamma r_2})} \le C \delta^{\alpha} \| u_\delta\|_{H^1(\Omega)}. 
\end{equation}
\medskip
In what follows, we assume that $k=1$ for notational ease.
Define $U_{2, \delta}$ in $B_{r_3} \setminus B_{\gamma r_2}$ as follows 
\begin{equation*}
U_{2, \delta} = u_{1, \delta} - u_{2, \delta} + V_{1, \delta}. 
\end{equation*}
Then  
\begin{equation}\label{eq-U2-Thm1} 
\Delta U_{2, \delta} +  U_{2, \delta} = 0 \quad  \mbox{ in } B_{r_3} \setminus B_{\gamma r_2},
\end{equation}
and
\begin{equation}\label{p-eq-U2}
U_{2, \delta} = 0 \quad\mbox{and }  \quad  \partial_r U_{2, \delta} = -\frac{\idelta}{1 - \idelta} \partial_r u_{1, \delta} + \partial_r V_{1, \delta} \quad \mbox{ on } \pB_{r_3}. 
\end{equation}

\medskip
\noindent \underline{Case 1:} $d=2$.  As in \cite{Ng-Superlensing}, define
\begin{equation*}
\hat J_n(r) = 2^n n! J_n(r) \quad \mbox{ and } \quad \hat Y_n(r) = \frac{\pi i}{2^{n} (n-1)!} Y_n(r), 
\end{equation*}
where $J_n$ and $Y_n$ are the Bessel and Neumann functions of order $n$. It follows from \cite[(3.57) and (3.58)]{ColtonKressInverse} that 
\begin{equation}\label{bh1}
\hat J_n (t)  = t^{n}\big[1 + o(1) \big]
\end{equation}
and
\begin{equation}\label{bh2}
\hat Y_n (t)  = t^{-n} \big[1 + o(1) \big], 
\end{equation}
as $n \to + \infty$.  
\medskip
From \eqref{eq-U2-Thm1} one can represent $U_{2, \delta}$ as follows
\begin{equation}\label{re-u-k}
U_{2, \delta} =  a_0 \hat J_0(|x|) + b_0 \hat Y_0(|x|) +   \sum_{n=1}^\infty \sum_{\pm} \big[a_{n, \pm} \hat J_{n}(|x|) + b_{n, \pm} \hat Y_{n}(|x|) \big] e^{\pm i n \theta} \quad \mbox{ in } B_{r_3} \setminus B_{\gamma r_2},
\end{equation}
for $a_0, b_0,  a_{n, \pm}, b_{n, \pm} \in \mC$ ($n \ge 1$). Assume that 
\begin{equation*}
\partial_r U_{2, \delta} = c_0 + \sum_{n \ge 1} \sum_{\pm} c_{n, \pm} e^{\pm i n \theta} \quad \mbox{ on } \partial B_{r_3}. 
\end{equation*}
Then, by \eqref{def-V1}, \eqref{V1}, and \eqref{p-eq-U2},
\begin{equation}\label{pro-c-n}
|c_0|^2 +  \sum_{n \ge 1} \sum_{\pm} n^{-1} |c_{n , \pm}|^2 \sim \|\partial_r U_{2, \delta} \|_{H^{-1/2}(\partial B_{r_3})}^2 \le C \delta^2 \| u \|_{H^{1}(\Omega)}^2. 
\end{equation}
Using \eqref{p-eq-U2} again, we have
\begin{equation*}
\left\{\begin{array}{c}
a_{n, \pm} \hat J_n (r_3) + b_{n, \pm} \hat Y_n (r_3) = 0 \\[6pt]
a_{n, \pm} \hat J_n' (r_3) + b_{n, \pm} \hat Y_n' (r_3) = c_{n, \pm}
\end{array}\right. \quad \mbox{ for } n \ge 0. 
\end{equation*}
Here we denote $a_{0, \pm} = a_0/2$; $b_{0, \pm} = b_0/2$; and $c_{0, \pm} = c_0/ 2$.  It follows that 
\begin{equation}\label{efcd-k}\left\{
\begin{array}{l}
\dsp a_{n, \pm}  = c_{n, \pm} AC_{n}, \\[6pt]
\dsp b_{n, \pm} = c_{n, \pm} BC_{n},   
\end{array} \right. \quad \mbox{ for } n \ge 0. 
\end{equation}
where
\begin{equation*}
AC_{n} = -  \frac{\hat Y_n}{\hat J_n \hat Y_n' - \hat J_n' \hat Y_n}(r_3)
\quad
\mbox{ and } 
\quad 
BC_{n} = - \frac{\hat J_n}{\hat Y_n \hat J_n' - \hat Y_n' \hat J_n}(r_3). 
\end{equation*}
Using \eqref{bh1} and \eqref{bh2}, we derive that 
\begin{equation*}
AC_{n, } = - \frac{1}{2n} r_3^{1-n} \Big( 1 + o(1) \Big) \quad \mbox{ and } BC_{n} = \frac{1}{2n} r_3^{1 + n} \Big( 1 + o(1) \Big). 
\end{equation*}
We now make use the removing of localized singularity technique introduced in \cite{Ng-Superlensing, Ng-Negative-Cloaking}. Set
\begin{equation*}
\hat u_\delta (x)=  \sum_{n=1}^\infty \sum_{\pm} b_{n, \pm} \hat Y_{n}(|x|) e^{\pm i n \theta} \quad \mbox{ in } B_{r_3} \setminus B_{\gamma r_2},
\end{equation*}
We claim that, for $\gamma r_2 \le r \le r_3$,  
\begin{equation}\label{U2-hu}
\|U_{2, \delta}  - \hat u_\delta \|_{H^{1/2}(\pB_{r})} + \|\partial_r U_{2, \delta}  - \partial_r  \hat u_\delta \|_{H^{-1/2}(\pB_{r})} \le C \delta \|u_\delta \|_{H^1(\Omega)}. 
\end{equation}
Indeed, for $\gamma r_2 \le r \le r_3$,  
\begin{align*}
\|U_{2, \delta}  - \hat u_\delta \|_{H^{1/2}(\pB_{r})}^2  = & \|\sum_{n \ge 0} \sum_{\pm} a_{n, \pm} \hat J_n(|x|) e^{i n \theta} \|_{H^{1/2}(\pB_{r})}^2 
 \sim   \sum_{n \ge 0} \sum_{\pm} (n+1)|a_{n, \pm}|^2 |\hat J_n(|x|)|^2 \\[6pt]
\sim  &  \sum_{n \ge 0} \sum_{\pm} (n+1)|c_{n, \pm} AC_n|^2  |\hat J_n(|x|)|^2  \le   C \sum_{n \ge 0} \sum_{\pm} (n+1)^{-1}| c_{n, \pm}|^2 (r/ r_3)^{2n}. 
\end{align*}
It follows from \eqref{pro-c-n} that 
\begin{equation*}
\|U_{2, \delta}  - \hat u_\delta \|_{H^{1/2}(\pB_{r})}  \le C \delta \| u_\delta \|_{H^1(\Omega)}, 
\end{equation*}
for $\gamma r_2 \le r \le r_3$. 
Similarly, 
\begin{equation*}
\|\partial_r U_{2, \delta}  - \partial_r \hat u_\delta \|_{H^{-1/2}(\pB_{r})}  \le C \delta \| u_\delta \|_{H^1(\Omega)}, 
\end{equation*}
for $\gamma r_2 \le r \le r_3$. As a consequence of \eqref{V1} and \eqref{U2-hu}, we obtain 
for $\gamma r_2 \le r \le r_3$,  
\begin{equation}\label{part2}
\|u_{1, \delta} - u_{2, \delta} - \hat u_\delta \|_{H^{1/2}(\pB_{r})} + \|\partial_r u_{1, \delta} - \partial u_{2, \delta}  - \partial_r  \hat u_\delta \|_{H^{-1/2}(\pB_{r})} \le C \delta \|u_\delta \|_{H^1(\Omega)}. 
\end{equation}
Define
\begin{equation*}
U_\delta = \left\{\begin{array}{cl} u_\delta & \mbox{ in } \Omega \setminus B_{r_3}, \\[6pt]
u_\delta - \hat u_\delta & \mbox{ if } x \in B_{r_3} \setminus B_{\gamma r_2}, \\[6pt]
u_{2, \delta} &  \mbox{ if } x \in B_{\gamma r_2}. 
\end{array} \right.
\end{equation*}
We have 
\begin{equation*}
\Div(A \nabla U_\delta) + k^2 \Sigma U_{\delta} = f \mbox{ in } \Omega \setminus (\partial B_{r_3} \cup \partial B_{\gamma r_2}). 
\end{equation*}
On the other hand, from \eqref{part1} and \eqref{part2}, we obtain 
\begin{equation*}
\|[U_\delta] \|_{H^{1/2}(\partial B_{\gamma r_2})} + \|[A \nabla U_\delta \cdot \nu ] \|_{H^{-1/2}(\partial B_{\gamma r_2})} \le C \delta^{\alpha} \| u_\delta \|_{H^1(\Omega)} 
\end{equation*}
and 
\begin{equation*}
\|[U_\delta] \|_{H^{1/2}(\partial B_{r_3})} + \|[A \nabla U_\delta \cdot \nu] \|_{H^{-1/2}(\partial B_{r_3})} \le C \delta^{\alpha} \| u_\delta \|_{H^1(\Omega)}. 
\end{equation*}
Using \eqref{haha0}, we derive that 
\begin{equation*}
\| U_\delta\|_{H^1(\Omega \setminus (\partial B_{r_3} \cup \partial  B_{\gamma r_2}) )} \le  C \delta^{\alpha} \Big( \delta^{-1/2}\| U_\delta \|_{L^2(\Omega \setminus B_{r_3})}^{1/2} \|f \|_{L^2(\Omega)}^{1/2} + \| f\|_{L^2(\Omega)}\Big) + C \|f \|_{L^2(\Omega)}.
\end{equation*}
Since $\alpha > 1/2$, it follows that $U_\delta$ is bounded in $H^1\big(\Omega \setminus (\partial B_{r_3} \cup \partial  B_{\gamma r_2}) \big)$. 
Without loss of generality,  one may assume that $U_{\delta} \to U$ weakly in $H^1\big(\Omega \setminus (\partial B_{r_3} \cup \partial  B_{\gamma r_2}) \big)$  as $\delta \to 0$; moreover, $U \in H^1(\Omega)$ and 
\begin{equation*}
\Delta U + k^2  U= f \mbox{ in } \Omega \mbox{ and } U = 0 \mbox{ on } \partial \Omega.
\end{equation*}
Hence $U = u$. Since the limit is unique, we have the convergence for the family $(U_{\delta})$ as $\delta \to 0$. 

\medskip
\noindent \underline{Case 2:} $d=3$.   Define 
\begin{equation*}
\hat j_n(t) =1 \cdot 3 \cdots (2n + 1) j_n(t) \quad \mbox{ and } \quad  \hat y_n = -  \frac{y_n(t)}{1 \cdot 3 \cdots (2n-1)} ,  
\end{equation*}
where $j_n$ and $y_n$ are the spherical Bessel and Neumann functions of order n. 
Then, as $n$ large enough,  (see,  e.g., \cite[(2.37) and (2.38)]{ColtonKressInverse})  
\begin{equation}\label{jy-n}
\hat j_n(kr) = r^n \big(1 + O(1/n) \big) \quad \mbox{ and } \quad \hat y_n(kr) = r^{-n-1} \big(1 + O(1/n) \big). 
\end{equation}
Thus one can represent $U_{2, \delta}$ of the form
\begin{equation}\label{re-u-k1}
U_{2, \delta} =  \sum_{n=1}^\infty \sum_{-n }^n \big[a^n_m \hat j_{n}(|x|) + b^n_m \hat y_{n}(|x|) \big] Y^n_m(\hat x ) \quad \mbox{ in } B_{r_3} \setminus B_{r_0},
\end{equation}
for $ a^n_m, b^n_m \in \mC$ and $\hat x = x/ |x|$. The proof now follows similarly as in the case $d=2$. The details are left to the reader. 
\proofend

\begin{remark} \fontfamily{m} \selectfont  In the proof, we use essentially the fact $(A, \Sigma) = (I, 1)$ in $B_{r_3} \setminus B_{\gamma r_2}$ to use separation of variables in this region. In fact, this condition is not necessary by using the technique of separation of variables for a general structure in \cite{Ng-CALR}. 
\end{remark}

\begin{remark} \label{rem-construction} \fontfamily{m} \selectfont 
The construction of the cloak given by \eqref{defA} is not restricted to the Kelvin transforms $F$ (and $G$). In fact, 
one can extend this construction to a general class of reflections considered in \cite{Ng-Complementary}.
\end{remark}

\begin{remark}  \fontfamily{m} \selectfont 
The condition $(F_*A, F_*\Sigma) = (A, \Sigma)$ in $B_{r_3} \setminus B_{r_2}$ is necessary to ensure that cloaking can be achieved and the localized resonance might take place see \cite{Ng-Helmholtz-WP} (see also \cite{AnneSophie-Ciarlet1} for related results).
\end{remark}

\begin{remark} \fontfamily{m} \selectfont  
Cloaking can also be achieved via schemes generated by changes of variables \cite{GreenleafLassasUhlmann,Leonhardt,PendrySchurigSmith}. 
Resonance might also appear in this context but  for  specific frequencies see \cite{KohnOnofreiVogeliusWeinstein, NguyenPerfectCloaking}.  It is shown in \cite{NguyenPerfectCloaking} that in the resonance case cloaking might not be achieved and the field inside the cloaked region can depend on the field outside.  Cloaking can also be achieved in the time regime via change of variables \cite{NguyenVogelius3, NguyenVogelius4}.  
\end{remark}

\section{Illusion optics using complementary media}\label{sect-illusion}
We next discuss briefly how to obtain  illusion optics in the spirit of  Lai et al. in \cite{LaiNgChenHanXiaoZhangChanIllusion}. The scheme used here is a combination of the ones  used for cloaking and superlensing in \cite{Ng-Negative-Cloaking, Ng-Superlensing} and is slightly different from \cite{LaiNgChenHanXiaoZhangChanIllusion}. More precisely, set 
\begin{equation*}
m = r_3^2/ r_2^2. 
\end{equation*}
Let $a_c \in [L^\infty(B_{r_2/m})]^{d \times d}$ be elliptic  and $\sigma_c \in L^{\infty}(B_{r_2^2/r_3^2}, \mC)$ with $\Im(\sigma_c) \geq 0$. Define
\begin{equation}\label{defA1}
A_1, \Sigma_1 = \left\{ \begin{array}{cl}  A, \Sigma & \mbox{ in } \Omega \setminus B_{r_2/ m}, \\[6pt]
a_c,  \sigma_c  & \mbox{ in } B_{r_2/ m}, 
\end{array} \right. 
\end{equation} 
and 
\begin{equation}\label{defA1hat}
\hat A_1, \hat \Sigma_1 = \left\{ \begin{array}{cl}  I, 1 \mbox{ in } & \Omega \setminus B_{r_2}, \\[6pt]
(r_3/ r_2)^{2-d} a_c(x/ m),  (r_3/ r_2)^{-d}\sigma_c(x/m)  & \mbox{ in } B_{r_2}. 
\end{array} \right. 
\end{equation} 
Recall that $(A, \Sigma)$ is defined in \eqref{defA}. 
We assume that the following equation has only zero solution in $H^1_0(\Omega)$:  
		\begin{equation}\label{P2}
			\Div (A_1\nabla v) + k^2 \Sigma_1 v = 0 \mbox{ in } \Omega. 
		\end{equation} 

We obtain the following result on illusion optics: 
\begin{theorem}\label{thm2} Let  $d=2, \, 3$, $f \in L^{2}(\Omega)$ with $\supp f \subset \Omega \setminus  B_{r_{3}}$ and let  $u$ and  $u_\delta$ in $ H^1_{0}(\Omega)$ be respectively the unique solution of 
\begin{equation*}
\dive(s_\delta A_1 \nabla u_\delta) + k^2 s_0 \Sigma_1 u_\delta = f \quad \mbox{in } \Omega,
\end{equation*}
and 
\begin{equation*}
\dive(\hat A_1 \nabla u) + k^2 \hat \Sigma_1 u= f \quad \mbox{in } \Omega. 
\end{equation*}
There exists $\gamma_0> 1$, depending {\bf only} on $\Lambda$ and the Lipschitz constant of $\hat a$  such that if $1 < \gamma < \gamma_0 $ then   
\begin{equation}\label{key-point-Sup}
u_{\delta} \to u \mbox{ weakly in } H^{1} (\Omega \setminus B_{r_3}) \mbox{ as } \delta \to 0. 
\end{equation}
\end{theorem}

For an observer outside $B_{r_3}$, the medium in $B_{r_3}$ looks like $(\hat A_1, \hat \Sigma_1)$: one has illusion optics. 

\medskip 

\noindent{\bf Proof.} The proof is similar to the one  of Theorem~\ref{thm1}. Note that in the proof of Theorem~\ref{thm1}, we do not use the information of the medium inside $B_{r_2/m}$. The details are left to the reader. \proofend




\end{document}